\documentclass[preprint,12pt]{elsarticle}

\usepackage{graphicx}
\usepackage{amsmath,amssymb,latexsym}
\usepackage[colorlinks, plainpages]{hyperref}
\usepackage[english,francais]{babel}

\newtheorem{theorem}{Theorem}[section]
\newtheorem{lemma}[theorem]{Lemma}

\newtheorem{corollary}[theorem]{Corollary}
\newtheorem{definition}[theorem]{Definition}

\newcommand{\be}{\begin{equation}}
\newcommand{\ee}{\end{equation}}
\newcommand{\beq}{\begin{eqnarray}}
\newcommand{\eeq}{\end{eqnarray}}

\begin{document}

\begin{frontmatter}

\title{ On a Relation between Spectral Theory of Lens Spaces and Ehrhart Theory}

\author{H. Mohades and B. Honari \fnref{}}

\begin{abstract}
In this article Ehrhart quasi-polynomials of simplices are employed to determine isospectral lens spaces in terms of a
finite set of numbers. Using the natural lattice associated with a lens space the associated toric variety of a lens space is introduced. It is proved that if two lens
spaces are isospectral then the dimension of global sections of powers of a natural line
bundle on these two toric varieties are equal and they have the same general intersection number. Also, harmonic polynomial representation of the group
SO($n$) is used to provide a direct proof for a theorem of Lauret, Miatello and Rossetti on isospectrality of lens spaces.
\end{abstract}


\end{frontmatter}



\section{Introduction}
\begin{par}
Laplace-Beltrami operator is a natural second order
elliptic operator on a Riemannian manifold defined as \textbf{div} $\circ$ \textbf{grad}. It is well known that on a closed manifold, this operator has discrete
positive eigenvalues with finite multiplicities \cite{Ca}.Two
Riemannian manifolds are isospectral if their Laplace-Beltrami
operator have the same spectrum, considering multiplicities. A
fundamental question by Mark Kac asks whether it is possible to find two
nonisometric isospectral manifolds. The first negative  answer to
this question was provided by Milnor's 16 dimensional tori that is a
geometric realization of self-dual lattices with the same theta
functions \cite{Mi}. Nowadays there are many methods to construct
such manifolds, including Vigneras' arithmetic method which benefits from
quaternion algebras \cite{V}, Sunada's triple method and its
generalizations based on free action of a triple $(G,H_1,H_2)$
consisting of a group G and two of its almost conjugate subgroups
that provide an equivalent linear representation  for regular
representations which are not permutationally equivalent \cite{Su,GW,G2},
Gordon's torus bundle method \cite{Bu} etc.
There are some examples of isospectral manifolds that cannot be derived from the above general methods
including isospectral lens spaces introduced by Ikeda \cite{I1,I2,I4,IT,IY}.

Necessity of calculation of analytic torsion led Ray to present explicitly the
spectrum of lens spaces. That was actually a continuation of works of Calabi, Gallot and Meyer
on spectrum of spheres \cite{R}.
Ikeda used representation theory to find the spectrum of lens spaces. He also
used his subtle method to study the spectrum of Hodge-Laplace operator
on differential forms \cite{I1,IY}. Representation theory methods has been applied by Ikeda, Gordon, Gornet, Lauret, McGowan, Miatello,
Rosetti and others to construct examples
of spherical orbifolds and space forms that clarify
differences between many types of isospectrality regarding the spectrum of
elliptic differential operators on these spaces \cite{GM,G2,LMR1,LMR2}. Recently Lauret, Miatello, and Rossetti extended representation theoretic
methods to find conditions on the natural lattices $\Gamma$
associated with isospectral lens spaces \cite{LMR1}. These conditions completely determine isospectral lens spaces.
Lauret also have recently applied the Ehrhart theory in more general context to drive a geometric characterization of isospectral orbifolds with the cyclic fundamental group \cite{L4}.
On the other hand, Ehrhart proved that the number of points of the  lattice $\mathbb{Z}^n$  which lay inside an integral multiple of a rational polygon can be obtained by a rational function \cite{BV}. In the case of a rational simplex, Macdonald, Stanley and others gave an explicit formula for this
rational function, so they provide an explicit formula for the number of lattice points
inside simple polytopes \cite{b34}.
Motivating from the theorem of Lauret, Miatello, and Rossetti, we naturally associate a simplex with a lens space. Using this simplex,
we introduce the Ehrhart polynomial and  associated toric variety of a lens space and we consider the effect of isospectrality of lens spaces on these objects.
For the sake of completeness, a direct proof of Lauret, Miatello, and Rosetti's theorem is provided by using the theory of harmonic polynomials as a realization for representation theory of SO($n$).

Here an overview of the paper is given.
In section 2 preliminaries on lens spaces is presented. Ehrhart polynomial is introduced in section 3. In section 4 a proof of Lauret, Miatello, and Rosetti's theorem is presented. In section 5 methods of previous sections are used to find conditions for isospectrality of lens spaces. Also the toric variety associated with a lens space is introduced and the relation between the spectrum of a lens space and the number of global sections of the natural line bundle over the toric variety is considered.
\end{par}
\section{Preliminaries on lens spaces}
\subsection{Harmonic Homogeneous polynomials}
A multivariate polynomial $P$ on $n$ variables with coefficients in $ \mathbb{C}$
is homogeneous of degree $m$ if and only if $P(\lambda x_1,...,\lambda
x_n)=\lambda^m P( x_1,...,x_n)$ for every $\lambda\in \mathbb{C}$.
The dimension of the complex vector space generated by these functions
is equal to the number of different monomials of degree $m$ in $n$
variables, i.e. $\left(\begin{array}{cc}
             n+m-1 \\
              n-1 \\
            \end{array}\right)$.
We denote  the $\mathbb{C}$-linear
space of homogeneous polynomials of degree $m$ on $\mathbb{R}^n$ with
coefficients in $\mathbb{C}$ and its subspace of harmonic homogeneous
functions by $P_m(\mathbb{R}^n)$ and $H_m(\mathbb{R}^n)$ respectively.  Let $\Delta$ = \textbf{div} $\circ$ \textbf{grad} be the Laplacian on the
Euclidean space $\mathbb{R}^n$ and  let $T:P_m(\mathbb{R}^n)\longrightarrow
P_m(\mathbb{R}^n)$ be a linear operator defined by $T(f)=r^2\Delta f$. Since Laplacian is an onto function from
$P_m (\mathbb{R}^n)$ to $P_{m-2} (\mathbb{R}^n)$, the kernel
of $T$ is equal to $H_m(\mathbb{R}^n)$ and its image equals $r^2
P_{m-2}(\mathbb{R}^n)$ . The following equation, therefore, holds

\begin {equation}
dimH_m=dimP_{m}-dimP_{m-2} = \left(\begin{array}{cc}
             n+m-1 \\
              m-1 \\
            \end{array}\right)-\left(\begin{array}{cc}
             n+m-3 \\
              m-3 \\
            \end{array}\right)
\end {equation}
\subsection{Lens spaces}
Let $q$ be a positive integer, and let $p_1,...,p_m$ be integers that are prime to $q$.
Let

\begin {equation}
R(\theta)=\left(
            \begin{array}{cc}
              cos\theta & -sin\theta \\
             sin\theta & cos\theta \\
            \end{array}
          \right)\sim e^{i\theta}
\end {equation}
 and

\begin {equation}
g=R(p_1/q)\oplus \cdots \oplus R(p_m/q).
\end {equation}
Suppose that $G \subset O(2m)$ is the finite cyclic group generated by $g$.
If $G$ as a group of isometries acts freely on $S^{2m
-1}$, then the manifold $L = S^{2m-1}/G$, denoted by $\L (q,p_1,...,p_m)$, is called a lens space.
Since the only finite groups with free action on even dimensional spheres are $ \mathbb{Z}_2$ and the
trivial group, we leave these spheres out.
 \subsection {Spherical harmonic functions over lens spaces}
We can extend a function $f(\theta)$ on $S^{n-1}$ to a function
$g(\theta)=f(\theta/|\theta|) $ on $\mathbb{R}^n-\{0\}$ and use the ordinary
Laplacian on $\mathbb{R}^n$ to define the operator $\Delta_{S^{n-1}}$ as
$\Delta_{S^{n-1}}f = \Delta g $. This operator equals \textbf{div} $\circ$
\textbf{grad} on $S^{n-1}$ with the induced inner product from $\mathbb{R}^n$.
Now we compute $\Delta$ in polar coordinates.\\
Let $r=|x|$ and $\theta=x/|x|$, then
\begin {equation} \label{Eqn:MMM}
\hspace{2 cm} \Delta( f(r)g(\theta)))=(\Delta f(r))g(\theta)+ f(r)\Delta(g(\theta)).
\end {equation}
We also have
\begin {equation}
 \frac{\partial r}{\partial x^i}=x^i/r \hspace{0.3cm} ,\hspace{0.3cm}  \frac{\partial^2 r}
{\partial(x^i)^2}=\frac{1}{r}-\frac{(x^i)^2}{r^3}
\end {equation}
 and

\begin {equation}
\Delta f(r) = \sum_\emph{i=0}^n \frac{\partial}{\partial x^i}
\frac{x^i}{r}    \frac{\partial f}{\partial r}=\sum_\emph{i=0}^n
\frac{\partial^2f}{\partial r^2} \frac{(x^i)^2}{r^2} +
\sum_\emph{i=0}^n ( \frac{1}{r }-\frac{(x^i)^2}{r^3})
\frac{\partial f}{\partial r}.
\end {equation}
Therefore
\begin {equation}
\Delta f = r^{1-n} \frac{\partial}{\partial r}\left(r^{n-1}\frac{\partial f}{\partial r}\right) + r^{-2} \Delta_{S^{n-1}}f.
\end {equation}
For $f(r) = r^k $ we obtain
\begin {equation}
\Delta (r^k g(\theta)) = r^{-2+k}(\Delta_{S^{n-1}} g + (k(k-1)+n k)g).
\end {equation}
If $  r^k g(\theta)$ is a
harmonic function on $\mathbb{R}^n-\{0\}$, then $g(\theta)$ is an
eigenfunction of $-\Delta_{S^{n-1}}$ with the eigenvalue $k(k+n-1)$.
Zero is a removable singularity for these harmonic functions and
therefore such harmonic functions are polynomials of degree $m$.

\section {Lattices and Ehrhart theory}
The well-known definitions and facts of this section will be used in subsection 4.2 and in section 5.
In this paper a lattice $L$ is considered as a subgroup of the group $\mathbb{Z}^n$.
Such a lattice is of full rank if $L\otimes \mathbb{R}=\mathbb{R}^n$.
Let $\{v_1,...,v_n\}$ be a basis for  $\textit L$.
The matrix whose columns are $v_1,...,v_n$ is called a generating matrix of the lattice $L$.
Two matrices $A$ and $B$ are generating matrices of the same lattice iff there is a unimodular matrix $U$ such that $A=UB$.
An essential parallelepiped is the parallelepiped $\{\sum_{i=1}^n a_i v_i | 0 \leq a_i \leq 1,i=1, \cdots,n\}$.
Define the lattice $L^*=\{x\in\mathbb{R}^n  |\,x\cdot y  \in  \mathbb{Z} , \forall y\in L\}$ to be the dual lattice  of the lattice $L$.
If $A$ is a generating matrix of $L$ then $(A^{-1})^T$ is a generating matrix for $L^*$.
A convex polytope in $\mathbb{R}^n$ is an integral (rational) polytope provided that all of its vertices are in
the lattice $\mathbb{Z}^n$ ($\mathbb{Q}^n $).\\
Let $H_n$ be the polytope  $\{x\in \mathbb{R}^n| |x|_{l_1}\leq 1\}$  and let $P$ be the  rational polytope $H_nA^{-1}$.
Also let $\{u_1,...,u_n\}$ be a basis of $L^*$.
Define $l(u_i)$ to be the smallest integer that $l(u_i)u_i \in \mathbb{Z}^n$ and set $\gamma_i=(l(u_i)u_i,l(u_i)) \in \mathbb{Z}^{n+1}$.
Define $\Delta_A=\{\sum_{i=1}^n a_i\gamma_i\hspace{0.1cm}|\hspace{0.1cm} 0\leq a_i  <1\}$  and
let $\Delta_A\cap \mathbb{Z}^{n+1}=\{(q_i,r_i)| q_i\in \mathbb{Z}^n, r_i\in \mathbb{Z}, 0\leq i< s\}$ where
$s$ is the greatest common divisor of all $n$ by $n$ minors of the matrix whose rows are $\gamma_i$.
Let $I(A,k)$ be the cardinal of the set $kH_nA^{-1}\cap \mathbb{Z}^n$ and set $J(A,x)=\sum_{k=1}^\infty I(A,k)x^k$.
A theorem of Stanley \cite{S34} asserts that $J(A,x)=(\sum_{i=0}^{s-1}x^{r_i})\prod_{j=0}^n(1-x^{l(u_j)})^{-1}$.
As a result $I(A,x)=\sum_{i=0}^n c_i(x)x^i$ where $c_i$ is a periodic function.
In the case that $P$ is a convex integral polytope, Ehrhart proved that $I(A,x)$ is a polynomial of $x$.
The leading term of Ehrhart polynomial is equal to $Vol(P)$ and the second coefficient is equal to
$\frac{1}{2}\sum_{\textbf{f}\in \textbf{F}} Vol(\textbf{f})$, when $\textbf{F}$ is the set of closed facets of the polytope $P$.

\section{Isospectrality of lens spaces}
\subsection{Eigenfunctions of lens spaces} Let $\pi:S^{2m-1}\rightarrow S^{2m-1}/G$ be the
natural local isometry from the sphere to the lens space $\mathfrak{L}=S^{2m-1}/G $.
Locally isometric spaces have the same local form of Laplace-Beltrami operator.
So if $f$  is an eigenfunction of Laplace-Beltrami operator on $S^{2m-1}/G$ associated with the eigenvalue $\lambda$,
then $ f \circ \pi$ is an eigenfunction of $-\Delta_{S^{2m-1}}$  with the same eigenvalue $\lambda$.
Multiplicity of an eigenvalue $\lambda$ is equal to the maximum number of independent eigenfunctions associated with the eigenvalue $\lambda$.
For the Laplace-Beltrami operator on a lens space this multiplicity is the same as the dimension of the space generated by those
$\lambda$-eigenfunctions of $-\Delta_{S^{2m-1}}$ which are invariant under the action of $G$ on $S^{2m-1}$.
So this is nothing but the dimension of the vector space of homogeneous harmonic polynomials on $\mathbb{R}^{2m}$
that are invariant under the action of $G$ on $S^{2m-1}$.

Note that  $G_1 \subseteq G$ implies  spec$({\mathbb{S}^n}/{G} )\subseteq$ spec$({\mathbb{S}^n}/{G_1})$.
In Particular, spec$({\mathbb{S}^n}/{G})$ $\subseteq$ spec $(\mathbb{S}^n )$.

\subsection{The lattice associated with a lens space}
The main part of this subsection is devoted to finding an elementary proof for the theorem of
Lauret, Miatello and Rossetti \ref{aq} and defining finitely many invariants which determine the isospectrality of two lens spaces.
The important fact that the Laplace-Beltrami operator on a manifold commutes with the isometries of
the manifold  and the main lemma \ref{uuu} are used to decompose the vector space of $\lambda$-eigenfunctions of the Laplace-Beltrami operator.
Let $z_j=x_j+iy_j$ and $z_{n+j}=\bar{z_j}$, $ j=1,\ldots ,n$.
In this case the monomials $\prod_{j=1}^n z_j^{\alpha_j }z_{n+j}^{\alpha_{n+j}} $, $\Sigma_{j=1}^{2n}\alpha_j=m$, make a basis for $P_m$.
This monomials can also be written as
\begin {equation} \label{kk} \prod_{j=1}^n|z_j |^{2\alpha_{\sigma(j)}} \prod_{j=1}^n z_{\sigma(j)}^{|\alpha_{n+j}-\alpha_j| }
\end {equation} where $\sigma(j)=n+j$
if $\alpha_{n+j}-\alpha_j \geq 0 $ and $\sigma(j)=j$ otherwise.
The action of an element  $\gamma\in G$ changes these monomials by a factor $e^{\frac{2\pi i}{q}\sum_{j=1}^n(\alpha_j-\alpha_{n+j}) p_j}$.
These monomials are linearly independent so a polynomial is invariant under the action of $G$  iff
all of its monomials are  invariant under this action.
This will happen iff
\begin{equation} \sum_{j=1}^n(\alpha_j-\alpha_{n+j}) p_j\equiv 0 \hspace{1cm}(\mathrm{mod}\hspace{0.1cm}q)
\end{equation}
 holds. \\
 \begin{definition}The lattice associated with a lens space $S^{2m-1}/G$ is defined to be the lattice $L=\{(a_1,...,a_n) \in \mathbb{Z}^n|\sum_{j=1}^n a_jp_j\equiv 0\hspace{0.2cm}(\mathrm{mod}\hspace{0.1cm}q)\}$\cite{LMR1}.
  \end{definition}
Note that there is a one to one correspondence between $L$ and the monomials of the form
$\prod_{j=1}^nz_{\sigma(j)}^{|\alpha_{n+j}-\alpha_j| }
$  that are invariant
under the action of $G$. Such correspondence can be figured out by mapping the element
$\alpha=(\alpha_{n+1}-\alpha_1
,\alpha_{n+2}-\alpha_2,...,\alpha_{2n}-\alpha_n )$ to the monomial $\prod_{j=1}^nz_{\sigma(j)}^{|\alpha_{n+j}-\alpha_j| }
$.\\
\begin {lemma}\label{uuu} Let each  $Q_l  $ be a  complex polynomial of $n$ variables. Let\\
$\sum_l Q_l(|z_1 |^2,...,|z_n|^2)\prod_{j=1}^nz_{\sigma(j)}^{|\alpha_{n+j}^l-\alpha_j^l| }$ be a harmonic function where the elements  $\prod_{j=1}^nz_{\sigma(j)}^{|\alpha_{n+j}^l-\alpha_j^l| }$ are distinct for different $l$'s. Then each  summand\\  $ Q_l(|z_1 |^2,...,|z_n
|^2)\prod_{j=1}^nz_{\sigma(j)}^{|\alpha_{n+j}^l-\alpha_j^l| } $ is harmonic.
\end {lemma}
Proof: The Laplacian commutes with the elements
 $$(R(\theta_1),R(\theta_2),...,R(\theta_n)) \in SO(2n).$$ So we have\\
 $$
\sum_j e^{i\sum_{l=1}^n({\theta_l{|\alpha_{n+l}^j-\alpha_l^j| })}} \Delta( Q_j(|z_1 |^2,...,|z_n
|^2)  \prod_{k=1}^n z_{\sigma(k)}^{|\alpha_{n+k}^j-\alpha_k^j| })=0.
$$ Now if for each $(\theta_1,...,\theta_n) \in \mathbb{R}^n$ the element $\sum_j e^{i\sum_{l=1}^n({\theta_l{|\alpha_{n+l}^j-\alpha_l^j| })}} \eta_j $  is equal to zero then $\eta_j$ must be zero. Thus the desired polynomials are harmonic.\\\\
Note that the homogeneous parts of harmonic polynomials are harmonic.
Now for the fixed element $\prod_{j=1}^n z_{\sigma (j)}^{|\alpha_{n+j} -\alpha_j|}$ we look for harmonic
polynomials of degree $m$ of the form
 \begin{equation} \label{ss} Q(|z_1 |^2,...,|z_n|^2) \prod_{j=1}^n z_{\sigma(j)}^{|\alpha_{n+j}-\alpha_j| }\hspace{0.2cm},
 \end {equation}
where $Q$ is a homogeneous polynomial of degree $ r$.
 The equation $2r+||\alpha||_1=m $ is satisfied by the degree $r$ and the element $\alpha$. We call  these kind of harmonic functions the harmonic polynomials associated with $\alpha$.
For the basis elements of $P_m$ let
 \begin {equation}  \label{dd00}
\prod_{j=1}^n |z_j |^{2\alpha_j }  z_{\sigma(j)}^{|\alpha_{n+j}-\alpha_j| } =p+R^2q
\end {equation} be the unique decomposition \cite{A} of the above monomial when $ p$ is a harmonic polynomial of degree $m$, $q$ is a homogeneous polynomial of degree $m-2$ and $R^2=|Z_1 |^2+...+|Z_n |^2$.
In this case, due to commutativity of the Laplacian with $\gamma \in G$, $\Delta( \gamma (p))$ is equal to zero and therefore $\gamma( p)$ will be a harmonic polynomial. The element $\gamma$ preserves the monomial, thus by the uniqueness of decomposition it preserves both $p$ and $q$.
By a similar argument, an arbitrary element $\eta$ of SO($n$), with the matrix form $(R(\theta_1),R(\theta_2),...,R(\theta_n))$, acts on two sides of the equation (\ref{dd00}) as follows
\begin {equation} \label{ss2}
e^{i \sum_{j=1}^n \theta_j |\alpha_{n+j}-\alpha_j| }\prod_{j=1}^n |z_j |^{2\alpha_{\sigma(j) }} z_{\sigma(j)}^{|\alpha_{n+j}-\alpha_j| }
=\eta p+R^2\eta q.
\end {equation}
Multiplying both sides of (\ref{dd00}) by  $e^{i \sum_{j=1}^n \theta_j |\alpha_{n+j}-\alpha_j| }$ we obtain
\begin {align} \label{ss1}
e^{i \sum_{j=1}^n \theta_j |\alpha_{n+j}-\alpha_j| }\prod_{j=1}^n |z_j |^{2\alpha_{\sigma(j) }} z_{\sigma(j)}^{|\alpha_{n+j}-\alpha_j| }
&=e^{i \sum_{j=1}^n \theta_j |\alpha_{n+j}-\alpha_j| } p \nonumber \\
& +R^2e^{i \sum_{j=1}^n \theta_j |\alpha_{n+j}-\alpha_j| } q.
\end {align}
Comparing (\ref{ss2}) and (\ref{ss1}), by uniqueness of decomposition we deduce$$\eta q=e^{i \sum_{j=1}^n \theta_j |\alpha_{n+j}-\alpha_j| } q.$$
Therefore each monomial of $q $ can be written in a form like (\ref{kk}).
It follows that\\
\begin {equation} \label{sre} q=h(|z_1 |^2,...,|z_n |^2)\prod_{j=1}^n z_{\sigma(j)}^{|\alpha_{n+j}-\alpha_j| } \end {equation}
where $h$ is a homogeneous polynomial of degree $r-1$. 
The following theorem calculates the dimension of these polynomials.
\begin {theorem}
The dimension of the vector space of  harmonic polynomials of degree $m$  associated with $\alpha$ is
\end {theorem}
$$\left(\begin{array}{cc}
             r+n-2 \\
              n-2 \\
            \end{array}\right).$$
\\
Proof: If we decompose each monomial of (\ref{ss}) and use the given format of (\ref{sre}) for $q$ we deduce that this dimension
is the same as the difference of dimensions of the spaces of homogeneous polynomials of degrees $r$ and $r-1$ on $n$ variables.\\\\
So for the lens space $\mathfrak{L}=S^{2n-1}/G$ the multiplicity of the eigenvalue $k(k+2n-1)-1)$ of the Laplace-Beltrami operator is equal to
\begin {equation} \label{tt43}
f_\mathfrak{L}(k)=\sum_{r=0}^{[k/2]} \left(\begin{array}{cc}
             r+n-2 \\
             n-2 \\
            \end{array}\right) N_L(k-2r)
\end {equation} where $N_L(s) $ is the number of elements $\alpha$ on the associated lattice $L$ in which $||\alpha||_1=s$.
Note that the  equality of $f_{\mathfrak{L}_1}$ and $f_{\mathfrak{L}_2}$ leads to the equality of $N_{L_1}$ and $N_{L_2}$ and vice versa \cite{LMR1}.
As a result we have the following theorem of Lauret, Miatello and Rossetti.
\begin {theorem} \label{aq}
Two lens spaces  $\mathfrak{L}_1=S^{2n-1}/G_1$ and $\mathfrak{L}_2=S^{2n-1}/G_2$ are isospectral
iff for the associated lattices $L_1$ and $L_2$, $N_{L_1}=N_{L_2}$.
\end {theorem}
For the lattice $L$, let $A$ be as in section 3. By  (\ref{tt43}) the multiplicity of the eigenvalue $k(k+n-1)$ is equal to\\
$$\sum_{r=0}^{[k/2]} \left(\begin{array}{cc}
             r+n-2 \\
              n-2 \\
            \end{array}\right) (I(A,k-2r)-I(A,k-2r-1)).$$
As a result of previuos theorem if the associated lattices of two lens spaces differ by an element of $SL(n,\mathbb{Z})$ then these lens spaces are isospectral.
\begin {theorem}\label{lenehr}
Two lens spaces with generating matrices of associated lattices $A$ and $B$  are isospectral iff their quasi-Ehrhart polynomials are equivalent;
i.e., $I(A,m)=I(B,m)$.
\end {theorem}
Proof: Trivially $I(A,0)=I(B,0)$. By Theorem \ref{aq}  we have $I(A,m)-I(A,m-1)=I(B,m)-I(B,m-1)$ for each $m\in \mathbb{N}$.
Therefore $I(A,m)-I(B,m)=0$ for each
            $m\in\mathbb{N}$. \\  

\begin {corollary} Referring to the notations of section 3, two lens spaces are isospectral iff the finitely many numbers $l(u_i)$
and $r_i$ are the same for the associated lattices.\end {corollary}
In order to determine a relation between lens spaces and toric varieties we need integral polytopes. Generally, $P=H_nA^{-1}$ is not an integral polytope. The following theorem 
leads us to modify our polytope $P$ to an appropriate one.
 \begin {lemma}\label{mul} Let $A$ and $B$ be the generating matrices of two Lattices.  Then, the  quasi-Ehrhart polynomials of
 convex rational polytopes $H_nA^{-1}$ and $H_nB^{-1}$ are equal iff for each $k \in \mathbb{N}$ the  quasi-Ehrhart polynomials of
 rational polytopes $k H_nA^{-1}$ and $k H_nB^{-1}$ are equal.
\end {lemma}
Proof: The cardinality of $xk H_nA^{-1}\cap \mathbb{Z}^n$  is equal to  $I(A,k x)$ which is equal to $I(B,k x)$.

\begin {theorem}\label{mul23}
 If two lens spaces are isospectral then\\
(i) The volume of
the fundamental domains of their associated lattices are equal.\\
(ii) The  boundary area of their associated convex polytopes are equal. 

\end {theorem}
Proof: $(i)$ Let $k\in \mathbb{N}$ be such that $kH_nA^{-1}$ and $kH_nB^{-1}$ are integral polytopes.
Then their Ehrhart polynomials are respectively equal  to $I(A,kx)$ and $I(B,kx)$.
These two polynomials are equal by Lemma \ref{mul}.
The leading terms of the given Ehrhart polynomials are volumes of convex polytopes $kH_nA^{-1}$ and $kH_nB^{-1}$.
These polytopes have volumes  $\frac {k^ndetA^{-1}}{n!}$ and $\frac {k^ndetB^{-1}}{n!}$ respectively.
Therefore $\det{A}=\det{B}$.
$(ii)$ The proof is the same as the one in the previous theorem, when we consider the second coefficient of the Ehrhart polynomial.
\section{ Toric varieties associated with lens spaces}
In this section we define the toric variety associated with a lens space and we consider the effect of isospectrality of two lens spaces on their toric varieties.
Our notations and definitions coincide with the ones in \cite{Co}.


Using the same notations as section 4, let $L$ be the associated lattice of the lens space $\mathfrak{L}$ with the generating matrix $A$ and let  $P$ denote  the polytope $(\det{A}) H_nA^{-1}$. 
An Integral polytope $P $ is very ample if for every vertex $m\in P$, the semigroup $N(P\cap \mathbb{Z}^n-m)$ generated by $\{m'-m|m' \in P\cap \mathbb{Z}^n\}$ is saturated in $\mathbb{Z}^n$; i.e., if for $k\in \mathbb{N}$, $kx\in N(P\cap \mathbb{Z}^n-m)$ then $x\in N(P\cap \mathbb{Z}^n-m)$.
\begin{theorem} Let $P $ be a full dimensional integral polytope of dimension $n \geq 2$. Then $nP$ is very ample \cite{Co}.
\end{theorem}

Note that very ampleness is precisely the property needed to define the toric variety of an integral polytope.

Let $P$ be a full n-dimensional very ample polytope 
and $\Sigma_P$  be the associated normal fan \cite{Fult}  of $P$. Let $\sigma ^*$ be the dual cone of $\sigma \in \Sigma_P$.
 Also, let $H_m=\{x\in \mathbb{Z}^n|\, x\,\cdot\, m=0\}$. If $m\in \sigma^*\cap \mathbb{Z}^n$, then for the face $\tau=\sigma\cap H_m$  of $\sigma$ we have $ Spec(\mathbb{C}[\tau^*\cap \mathbb{Z}^n])$ is isomorphic to the localization $ Spec(\mathbb{C}[\sigma^*\cap \mathbb{Z}^n])_{\chi^m}$, where $\chi^m:=\prod_{i=1}^n x_i^{m_i}$. Gluing the affine varieties $Spec(\mathbb{C}[\sigma^*\cap \mathbb{Z}^n])$  along the subvarieties $(Spec(\mathbb{C}[\sigma^*\cap \mathbb{Z}^n]))_{\chi^m}$,  we obtain an abstract variety $X_{\Sigma_P}$ which is a normal toric variety.
According to Lemma \ref{mul}, if two lens spaces are isospectral then the Ehrhart polynomials of the same multiple of their associated polytopes are equal. Motivating by these facts the following definition is given.

\begin {definition}A toric variety of a lens space is the toric variety $X_{\Sigma_P}$ constructed by the very ample polytope $nP$.\end{definition}
 $ \mathbb{C}^{\times n}$ is the maximal algebraic torus which acts freely on a dense subset of the toric variety. Let  $ \Sigma(1)$  be the set of 1-dimensional cones in $\Sigma_P$. The 1-dimensional cones $\rho \in \Sigma(1)$ correspond to $ (n-1)$-dimensional $ \mathbb{C}^{\times n}$-orbits in $X_{\Sigma_P}$ denoted by $O(\rho)$.  Let $D_\rho$  be the Zariski closure of the orbit $O(\rho)$, then $D_\rho$ is a $ \mathbb{C}^{\times n}$--invariant prime divisor of $X_\Sigma$.
Let $D=\sum_{\rho\in \Sigma(1)} a_\rho D_\rho$ be a torus invariant Cartier divisor on $X_{\Sigma_P}$.
For each $\rho \in \Sigma(1)$, we can choose $a_\rho$ naturally by real numbers used to define the $\rho$-facets of $P$ \cite{Co}.
So $D$ will be a natural divisor. Now for the sections of the sheaf $\emph{O}_ {X_{\Sigma_P}}(D)$  we have
$\Gamma (X_{\Sigma_P},\emph{O}_ {X_{\Sigma_P}}(D))=\bigoplus_{m\in P \cap \mathbb{Z}^n}(\mathbb{C}\bigotimes_\mathbb{Z}<\chi^m>)$, where $<\chi^m> $ is the $\mathbb{Z}$--module generated by $\chi^m$ \cite{Co}.

We know that on a normal variety $X$, the sheaf $O_X(D)$ is equal to the sheaf of sections of a line bundle $L$ which is unique up to isomorphism.
We call $L$, the natural line bundle on the toric variety $X_P$.
Let $H^0 (X_P,L^k )$ denote the sections of the $k$th power of $L$.
The above discussion leads us to the next theorem for toric varieties associated with lens spaces.
\begin {theorem} \label{78} Let $\mathfrak{L}_1$ and $\mathfrak{L}_2$ be two lens spaces and let $X_1$ and $X_2$ be the toric varieties associated with them.
If $\mathfrak{L}_1$ and $\mathfrak{L}_2$ are isospectral, then \\
(i) $dim H^0 (X_1,L^k )= dim H^0 (X_2,L^k )$ for each $k\in \mathbb{N}$.\\
(ii) $deg(X_1)=deg(X_2)$
\end {theorem}Proof: $(i)$ The dimension of global sections of $L^k$ is equal to the cardinality of $kP\cap \mathbb{Z}^n$. So the proof is a result of Theorem \ref{lenehr}. $(ii)$  From the Bernstein-Kouchnirenko theorem, we know that the $\frac{1}{n!}deg(X_i)$ are equal to the leading term of Hilbert polynomial of $X_i$ \cite{kav} and this leading term is equal to the volume of the associated polytope. So the theorem is a result of Theorem \ref{mul23}.\\\\
\textbf{Question.1.} Find a natural toric variety $X$ (integral polytope) associated with a lens space, in such a way that the isospectrality of two lens spaces be equivalent to the equality on  $dim H^0 (X ,L^k )$ for each $k\in \mathbb{N}$.\\
\textbf{Question.2.} Generalize the theorem \ref{78} when isospectrality being considered with respect to Hodge-Laplace operators.



\begin{thebibliography}{9}

%
%

\bibitem{A} S. Axler,  \textit{Harmonic Function Theory},  Springer, 2001.


\bibitem{b34}A. Barvinok,  \textit{ Computing the Ehrhart quasi-polynomial of a rational simplex}, Mathematics of Computation, 75(255), 1449-1466, (2006).
\bibitem{Bu}P. Buser,  \textit{Geometry and Spectra of Compact Riemann Surfaces}, Birkhauser, 2010.
\bibitem{Ca}I. Chavel, \textit{Eigenvalues in Riemannian Geometry}, Academic Press, 1984.
\bibitem{Co}D.  Cox, J. Little,   H. Schenck, \textit{Toric Varieties}, Vol. 124, American Mathematical Soc., 2011.

\bibitem{DG}E. Dryden,  V. Guillemin, and R. Sena-Dias,  \textit{Equivariant inverse spectral theory and toric orbifolds}, Advances in Mathematics, 231(3), 1271-1290, (2012).
\bibitem{BV} E. Ehrhart, \textit{ Polynomes arithmetiques et Methode des Polyedres en Combinatoire,} Internat.
Ser. Numer. Math, vol. 35, Birkhauser, Basel, (1977).
\bibitem{Fult}W. Fulton,  \textit{ Introduction to Toric Varieties}, Annals of Mathematics Studies 131, Princeton Univ.
Press, Princeton, NJ, 1993.

\bibitem{Z}C. Gordon, D. Webb,  S. Wolpert,  \textit{Isospectral plane domains and surfaces via Riemannian orbifolds},  Inventiones mathematicae, 110.1, 1-22, (1992).
\bibitem{GW}C. Gordon, D. Webb,  \textit{You can't hear the shape of a drum}. , American Scientist 84 (January-February), 46–55 (1996).
\bibitem{G2}C. Gordon,  \textit {Sunada's isospectrality technique: Two decades later,}  Spectral Analysis in Geometry and Number Theory,   484, 45-58, (2009).
 \bibitem {G1}C. Gordon,\textit{ Riemannian manifolds isospectral on functions but not on 1-forms}, J. Diff. Geom, 24, 79-96, (1986).
 \bibitem{GM}R. Gornet, J. McGowan, \textit{ Lens spaces, isospectral on forms but not on functions,} LMS Journal of Computation and Mathematics, Volume 9, 270-286, (2006).
\bibitem{IY}A. Ikeda,  Y. Yamamoto,  \textit{On the spectra of 3-dimensional lens spaces}, Osaka J. Math. 16, 447-469, (1979).
\bibitem{IT}A. Ikeda,  Y. Taniguchi,  \textit{ Spectra and eigenforms of the Laplacian on $S^n$ and $P^n(C)$}.    Osaka J. Math.Volume 15, Number 3, 515-546, (1978).
\bibitem{I1}A. Ikeda,  \textit{Riemannian manifolds p-isospectral but not p + 1-isospectral}, Geometry
of manifolds (Matsumoto, 1988), Perspect. Math. 8, 383-417, (1989).
\bibitem{I2}A. Ikeda,  \textit{On spherical space forms which are isospectral but not isometric}, J.
Math. Soc, Japan 35:3, 437-444, (1983).
\bibitem{I4}A. Ikeda, \textit{On lens spaces which are isospectral but not isometric},
 Ann. Sci. ´ Ecole
Norm. Sup. (4), 13:3, 303-315, (1980).
\bibitem{kav}K. Kaveh, A.  Khovanskii,  \textit{Algebraic equations and convex bodies}, In Perspectives in analysis, geometry, and topology, 263-282, Birkhauser Boston, (2012).
 \bibitem{L4}E.A. Lauret, \textit{ Spectra of orbifolds with cyclic fundamental groups,} Ann Glob Anal Geom, (2016).
 
\bibitem{LMR1}E.A. Lauret, R. J. Miatello, J. P. Rossetti,  \textit{Lens spaces isospectral on p-forms for every p,}
 arXiv, (2013).

\bibitem{LMR2}E.A. Lauret, R. J. Miatello, J. P. Rossetti,\textit{ Representation equivalence and p-spectrum of constant curvature space forms,} arXiv, (2012).
\bibitem{Macdonald}I.G. Macdonald, \textit{Polynomials associated with Finite cell-complexes}, Journal of the London Mathematical Society, 2(1), 181-192, (1971).
\bibitem{Mi}J. Milnor, \textit{Eigenvalues of the Laplace operator on certain manifolds}, Proc. Nat. Acad. Sci, USA, 51(4), 542, (1964).


\bibitem{R}D.B. Ray, \textit{Reidemeister torsion and the laplacian on lens spaces}, Advances in Mathematics 4, 109-126, (1970).



\bibitem{P} P. Sole, \textit{Counting lattice points in pyramids}, Bull. Sc. math, $2^e$ s\'{e}rie, 100, 149-173, (1976).
\bibitem{S34}R. P. Stanley, \textit{Decompositions of rational convex polytopes}, Ann. Discrete Math, v6, 333-342, (1980).
\bibitem{Su}T. Sunada, \textit{Riemannian coverings and isospectral manifolds}, Ann. math, (1985).
\bibitem{V}W. Vigneras,  \textit{The arithmetic of quaternion algebra}, preprint, 2006.
\bibitem{W}J.A. Wolf,  \textit{ Isospectrality for spherical space forms}, Result. Math, 40, 321-338, (2001).
\bibitem{Z}W. Ziller, \textit{Lie groups, representation theory and
symmetric spaces}, University of Pennsylvania, preprint, 2010.



\end{thebibliography}
\end{document}